# Numeration systems as dynamical systems – introduction

**Teturo Kamae**[1]

*Matsuyama University*

**Abstract:** A numeration system originally implies a digitization of real numbers, but in this paper it rather implies a compactification of real numbers as a result of the digitization.

By definition, a numeration system with $G$, where $G$ is a nontrivial closed multiplicative subgroup of $\mathbb{R}_+$, is a nontrivial compact metrizable space $\Omega$ admitting a continuous $(\lambda\omega + t)$-action of $(\lambda, t) \in G \times \mathbb{R}$ to $\omega \in \Omega$, such that the $(\omega + t)$-action is strictly ergodic with the unique invariant probability measure $\mu_\Omega$, which is the unique $G$-invariant probability measure attaining the topological entropy $|\log \lambda|$ of the transformation $\omega \mapsto \lambda\omega$ for any $\lambda \neq 1$.

We construct a class of numeration systems coming from weighted substitutions, which contains those coming from substitutions or $\beta$-expansions with algebraic $\beta$. It also contains those with $G = \mathbb{R}_+$.

We obtained an exact formula for the $\zeta$-function of the numeration systems coming from weighted substitutions and studied the properties. We found a lot of applications of the numeration systems to the $\beta$-expansions, Fractal geometry or the deterministic self-similar processes which are seen in [10].

This paper is based on [9] changing the way of presentation. The complete version of this paper is in [10].

## 1. Numeration systems

By a *numeration system*, we mean a compact metrizable space $\Omega$ with at least 2 elements as follows:

($\sharp$1) There exists a nontrivial closed multiplicative subgroup $G$ of $\mathbb{R}_+$ and a continuous action $\lambda\omega + t$ of $(\lambda, t) \in G \times \mathbb{R}$ to $\omega \in \Omega$ such that $\lambda'(\lambda\omega + t) + t' = \lambda'\lambda\omega + \lambda' t + t'$.

($\sharp$2) The $(\omega + t)$-action of $t \in \mathbb{R}$ to $\omega \in \Omega$ is strictly ergodic with the unique invariant probability measure $\mu_\Omega$ called the *equilibrium measure* on $\Omega$. Consequently, it is invariant under the $(\lambda\omega + t)$-action of $(\lambda, t) \in G \times \mathbb{R}$ to $\omega \in \Omega$ as well.

($\sharp$3) For any fixed $\lambda_0 \in G$, the transformation $\omega \mapsto \lambda_0\omega$ on $\Omega$ has the $|\log \lambda_0|$-topological entropy. For any probability measure $\nu$ on $\Omega$ other than $\mu_\Omega$ which is invariant under the $\lambda\omega$-action of $\lambda \in G$ to $\omega$, and $1 \neq \lambda_0 \in G$, it holds that

$$h_\nu(\lambda_0) < h_{\mu_\Omega}(\lambda_0) = |\log \lambda_0|.$$

The $(\omega + t)$-action of $t \in \mathbb{R}$ to $\omega \in \Omega$ is called the *additive* action or $\mathbb{R}$-action, while the $\lambda\omega$-action of $\lambda \in G$ to $\omega \in \Omega$ is called the *multiplicative* action or $G$-action.

Note that if $\Omega$ is a numeration system, then $\Omega$ is a connected space with the continuum cardinality. Also, note that the multiplicative group $G$ as above is either $\mathbb{R}_+$ or $\{\lambda^n;\ n \in \mathbb{Z}\}$ for some $\lambda > 1$. Moreover, the additive action is faithful, that is, $\omega + t = \omega$ implies $t = 0$ for any $\omega \in \Omega$ and $t \in \mathbb{R}$.

[1]Matsuyama University, 790-8578 Japan, e-mail: kamae@apost.plala.or.jp

This is because if there exist $\omega_1 \in \Omega$ and $t_1 \neq 0$ such that $\omega_1 + t_1 = \omega_1$, then take a sequence $\lambda_n$ in $G$ such that $\lambda_n \to 0$ and $\lambda_n \omega_1$ converges as $n \to \infty$. Let $\omega_\infty := \lim_{n\to\infty} \lambda_n \omega_1$. For any $t \in \mathbb{R}$, let $a_n$ be a sequence of integers such that $a_n \lambda_n t_1 \to t$ as $n \to \infty$. Then we have

$$\begin{aligned}\omega_\infty + t &= \lim_{n\to\infty} (\lambda_n \omega_1 + \lambda_n a_n t_1) \\ &= \lim_{n\to\infty} \lambda_n (\omega_1 + a_n t_1) = \lim_{n\to\infty} \lambda_n \omega_1 = \omega_\infty.\end{aligned}$$

Thus, $\omega_\infty$ becomes a fixed point of the $(\omega + t)$-action of $t \in \mathbb{R}$ to $\omega \in \Omega$. Since this action is minimal, we have $\Omega = \{\omega_\infty\}$, contradicting with that $\Omega$ has at least 2 elements.

An example of a numeration system is the set $\{0, 1\}^{\mathbb{Z}}$ with the product topology divided by the closed equivalence relation $\sim$ such that

$$(\dots, \alpha_{-2}, \alpha_{-1}; \alpha_0, \alpha_1, \alpha_2, \dots) \sim (\dots, \beta_{-2}, \beta_{-1}; \beta_0, \beta_1, \beta_2, \dots)$$

if and only if there exists $N \in \mathbb{Z} \cup \{\pm\infty\}$ satisfying that $\alpha_n = \beta_n$ $(\forall n > N)$, $\alpha_N = \beta_N + 1$ and $\alpha_n = 0$, $\beta_n = 1$ $(\forall n < N)$ or the same statement with $\alpha$ and $\beta$ exchanged. Let $\Omega(2) := \{0, 1\}^{\mathbb{Z}} / \sim$ and the equivalence class containing $(\dots, \alpha_{-2}, \alpha_{-1}; \alpha_0, \alpha_1, \alpha_2, \dots) \in \{0, 1\}^{\mathbb{Z}}$ is denoted by $\sum_{n=-\infty}^{\infty} \alpha_n 2^n \in \Omega(2)$. Then, $\Omega(2)$ is an additive topological group with the addition as follows:

$$\sum_{n=-\infty}^{\infty} \alpha_n 2^n + \sum_{n=-\infty}^{\infty} \beta_n 2^n = \sum_{n=-\infty}^{\infty} \gamma_n 2^n$$

if and only if there exists $(\dots, \eta_{-2}, \eta_{-1}; \eta_0, \eta_1, \eta_2, \dots) \in \{0, 1\}^{\mathbb{Z}}$ satisfying that

$$2\eta_{n+1} + \gamma_n = \alpha_n + \beta_n + \eta_n \quad (\forall n \in \mathbb{Z}).$$

This is isomorphic to the 2-adic *solenoidal group* which is by definition the projective limit of the projective system $\theta : \mathbb{R}/\mathbb{Z} \to \mathbb{R}/\mathbb{Z}$ with $\theta(\alpha) = 2\alpha$ $(\alpha \in \mathbb{R}/\mathbb{Z})$.

Moreover, $\mathbb{R}$ is imbedded in $\Omega(2)$ continuously as a dense additive subgroup in the way that a nonnegative real number $\alpha$ is identified with $\sum_{n=-\infty}^{\infty} \alpha_n 2^n$ such that $\alpha = \sum_{n=-\infty}^{N} \alpha_n 2^n$ and $\alpha_n = 0$ $(\forall n > N)$ for some $N \in \mathbb{Z}$, while a negative real number $-\alpha$ with $\alpha$ as above is identified with $\sum_{n=-\infty}^{\infty} (1 - \alpha_n) 2^n$. Then, $\mathbb{R}$ acts additively to $\Omega(2)$ by this addition. Furthermore, $G := \{2^k;\ k \in \mathbb{Z}\}$ acts multiplicatively to $\Omega(2)$ by

$$2^k \sum_{n=-\infty}^{\infty} \alpha_n 2^n = \sum_{n=-\infty}^{\infty} \alpha_{n-k} 2^n.$$

Thus, we have a group of actions on $\Omega(2)$ satisfying (♯1), (♯2)and (♯3) with $G := \{2^k;\ k \in \mathbb{Z}\}$ and the equilibrium measure $(1/2, 1/2)^{\mathbb{Z}}$.

**Theorem 1.1.** $\Omega(2)$ *is a numeration system with* $G = \{2^n;\ n \in \mathbb{Z}\}$.

We can express $\Omega(2)$ in the following different way. By a *partition* of the upper half plane $\mathbb{H} := \{z = x + iy;\ y > 0\}$, we mean a disjoint family of open sets such that the union of their closures coincides with $\mathbb{H}$. Let us consider the space $\Omega(2)'$ of partitions $\omega$ of $\mathbb{H}$ by open squares of the form $(x_1, x_2) \times (y_1, y_2)$ with $x_2 - x_1 = y_2 - y_1 = y_1$ and $y_1 \in G$ such that $(x_1, x_2) \times (y_1, y_2) \in \omega$ implies

$$\begin{aligned}&(x_1, (x_1 + x_2)/2) \times (y_1/2, y_1) \in \omega \quad \text{(type 0)} \quad \text{and} \\ &((x_1 + x_2)/2, x_2) \times (y_1/2, y_1) \in \omega \quad \text{(type 1)}.\end{aligned} \tag{1}$$

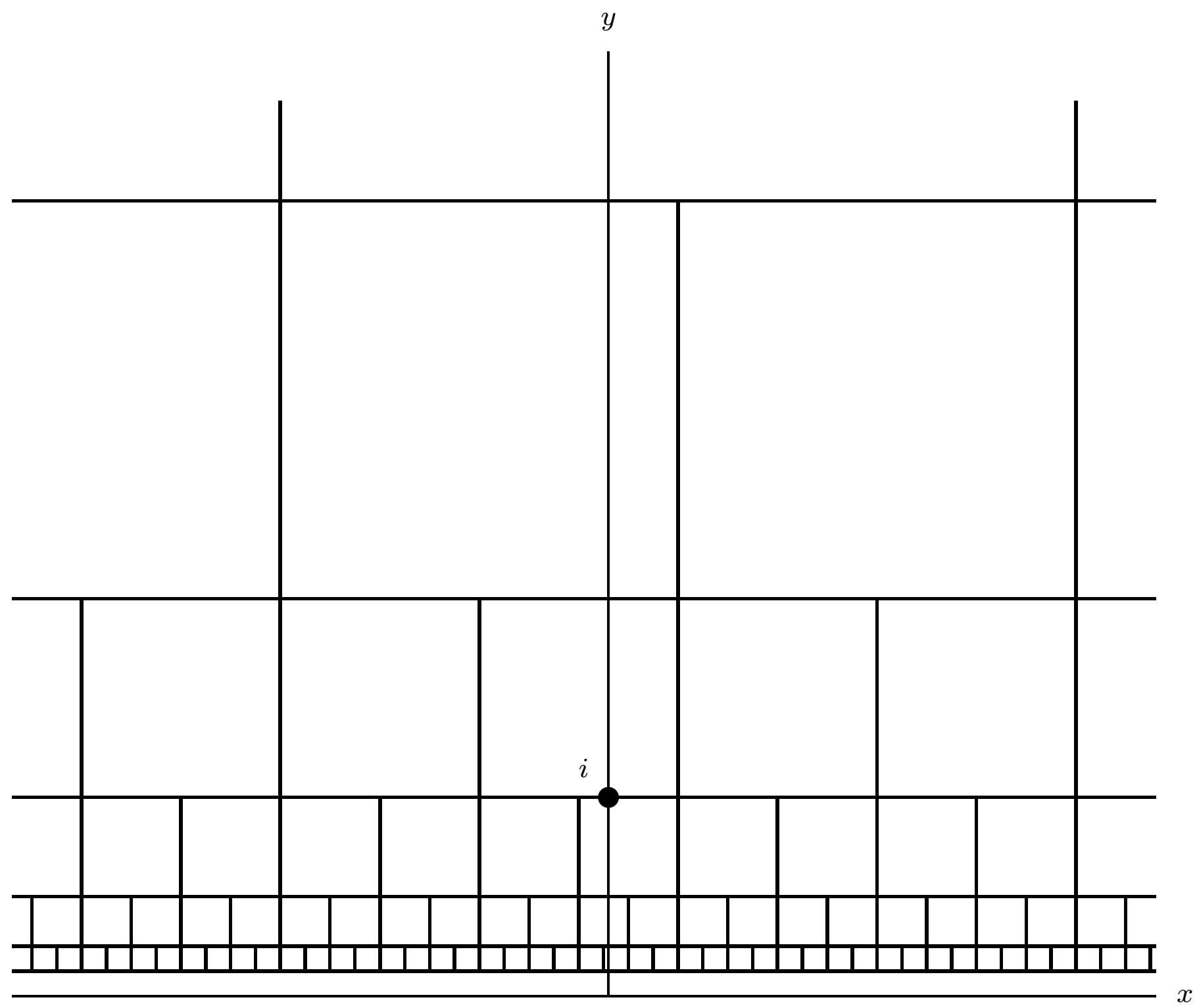

FIG 1. *The tiling corresponding to* $\cdots 01.101 \cdots$.

An example of $\omega \in \Omega(2)'$ is shown in Figure 1. For $\omega \in \Omega(2)'$, let $(\alpha_0, \alpha_1, \ldots)$ be the sequence of the types defined in (1) of the squares in $\omega$ intersecting with the half vertical line from $+0+i$ to $+0+i\infty$ and let $(\alpha_{-1}, \alpha_{-2}, \ldots)$ be the sequence of the types of the squares in $\omega$ intersecting with the line segment from $+0+i$ to $+0$. Then, $\omega$ is identified with $\sum_{n=-\infty}^{\infty} \alpha_n 2^n$. Note that replacing $+0$ by $-0$, we get $\sum_{n=-\infty}^{\infty} \beta_n 2^n$ such that $(\ldots, \alpha_{-2}, \alpha_{-1}; \alpha_0, \alpha_1, \ldots) \sim (\ldots, \beta_{-2}, \beta_{-1}; \beta_0, \beta_1, \ldots)$.

The topology on $\Omega(2)'$ is defined so that $\omega_n \in \Omega(2)'$ converges to $\omega \in \Omega(2)'$ as $n \to \infty$ if for every $R \in \omega$, there exist $R_n \in \omega_n$ such that $\lim_{n\to\infty} \rho(R, R_n) = 0$, where $\rho$ is the Hausdorff metric between sets $R,\ R' \subset \mathbb{H}$

$$\rho(R, R') := \max\{\sup_{z\in R} \inf_{z'\in R'} |z - z'|,\ \sup_{z'\in R'} \inf_{z\in R} |z - z'|\}. \tag{2}$$

For $\omega \in \Omega(2)'$, $t \in \mathbb{R}$ and $\lambda \in \{2^n;\ n \in \mathbb{R}\}$, $\omega + t \in \Omega(2)'$ and $\lambda\omega \in \Omega(2)'$ are defined as the partitions

$$\omega + t := \{(x_1 - t, x_2 - t) \times (y_1, y_2);\ (x_1, x_2) \times (y_1, y_2) \in \omega\}$$

and

$$\lambda\omega := \{(\lambda x_1, \lambda x_2) \times (\lambda y_1, \lambda y_2);\ (x_1, x_2) \times (y_1, y_2) \in \omega\}.$$

Let $\kappa : \Omega(2)' \to \Omega(2)$ be the identification mapping defined above. Then, $\kappa$ is a homeomorphism between $\Omega(2)'$ and $\Omega(2)$ such that $\kappa(\omega + t) = \kappa(\omega) + t$ and $\kappa(\lambda\omega) = \lambda\kappa(\omega)$ for any $\omega \in \Omega(2)'$, $t \in \mathbb{R}$ and $\lambda \in \{2^n;\ n \in \mathbb{Z}\}$. Thus, $\Omega(2)'$ is isomorphic to $\Omega(2)$ as a numeration system and will be identified with $\Omega(2)$.

We generalize this construction. Let $\mathbb{A}$ be a nonempty finite set. An element in $\mathbb{A}$ is called a *color*. An open rectangle $(x_1, x_2) \times (y_1, y_2)$ in $\mathbb{H}$ is called an *admissible tile* if

$$x_2 - x_1 = y_1 \tag{3}$$

is satisfied (see Figure 2). In another word, an admissible tile is a rectangle $(x_1, x_2) \times (y_1, y_2)$ in $\mathbb{H}$ such that the lower side has the hyperbolic length 1. Let $\mathcal{R}$ be the set of admissible tiles in $\mathbb{H}$.

A *colored tiling* $\omega$ is a subset of $\mathcal{R} \times \mathbb{A}$ such that
(1) $R \cap R' = \emptyset$ for any $(R, a)$ and $(R', a')$ in $\omega$ with $(R, a) \neq (R', a')$, and
(2) $\cup_{a \in \mathbb{A}} \cup_{(R,a) \in \omega} \overline{R} = \mathbb{H}$.

An element in $\mathcal{R} \times \mathbb{A}$ is called a *colored tile*. We denote

$$dom(\omega) := \{R;\ (R, a) \in \omega \text{ for some } a \in \mathbb{A}\}.$$

For $R \in dom(\omega)$, there exists a unique $a \in \mathbb{A}$ such that $(R, a) \in \omega$, which is denoted by $\omega(R)$ and is called the *color* of the tile $R$ (in $\omega$). Let $R = (x_1, x_2) \times (y_1, y_2)$. We call $y_2/y_1$ the *vertical size* of the tile $R$ which is denoted by $S(R)$.

Let $\Omega(\mathbb{A})$ be the set of colored tilings with colors in $\mathbb{A}$. A topology is introduced on $\Omega(\mathbb{A})$ so that a net $\{\omega_n\}_{n \in I} \subset \Omega(\mathbb{A})$ converges to $\omega \in \Omega(\mathbb{A})$ if for every $(R, a) \in \omega$, there exists $(R_n, a_n) \in \omega_n$ such that

$$a_n = a \text{ for any sufficiently large } n \in I \text{ and } \lim_{n \to \infty} \rho(R, R_n) = 0,$$

where $\rho$ is the Hausdorff metric defined in (2).

For an admissible tile $R := (x_1, x_2) \times (y_1, y_2)$, $t \in \mathbb{R}$ and $\lambda \in \mathbb{R}_+$, we denote

$$\begin{aligned} R + t &:= (x_1 + t, x_2 + t) \times (y_1, y_2) \\ \lambda R &:= (\lambda x_1, \lambda x_2) \times (\lambda y_1, \lambda y_2). \end{aligned}$$

Note that they are also admissible tiles.

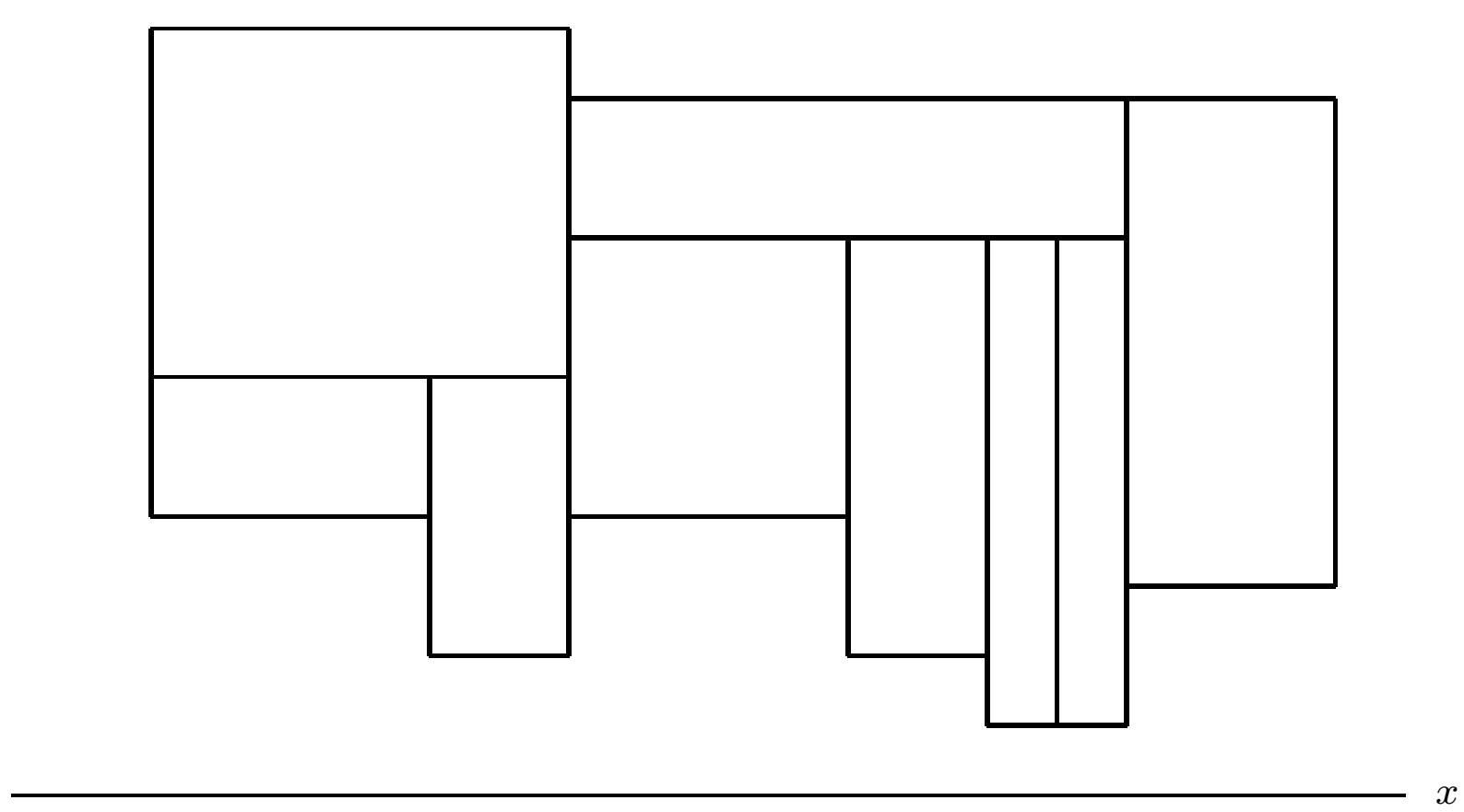

Fig 2. *Admissible tiles.*

For $\omega \in \Omega(\mathbb{A})$, $t \in \mathbb{R}$ and $\lambda \in \mathbb{R}_+$, we define $\omega + t \in \Omega(\mathbb{A})$ and $\lambda\omega \in \Omega(\mathbb{A})$ as follows:

$$\begin{aligned}\omega + t &= \{(R - t, a);\ (R, a) \in \omega\}\\ \lambda\omega &= \{(\lambda R, a);\ (R, a) \in \omega\}.\end{aligned}$$

Thus, we define a continuous group action $\lambda\omega + t$ of $(\lambda, t) \in \mathbb{R}_+ \times \mathbb{R}$ to $\omega \in \Omega(\mathbb{A})$. We construct compact metrizable subspaces of $\Omega(\mathbb{A})$ corresponding to weighted substitutions which are numeration systems. Though $\sharp\mathbb{A} \geq 2$ is assumed in [7], we consider the case $\sharp\mathbb{A} = 1$ as well.

## 2. Remarks on the notations

In this paper, the notations are changed in a large scale from the previous papers [7], [8] and [9] of the author. The main changes are as follows:

(1) Here, the colored tilings are defined on the upper half plane $\mathbb{H}$, not on $\mathbb{R}^2$ as in the previous papers. The multiplicative action here agree with the multiplication on $\mathbb{H}$, while it agree with the logarithmic version of the multiplication at one coodinate in the previous papers. Here, the tiles are open rectangles, not half open rectangles as in the previous papers.

(2) Here, we simplified the proof in [9] for the space of colored tilings coming from weighted substitutions to be numeration systems by omitting the arguments on the topological entropy.

(3) The roles of $x$-axis and $y$-axis for colored tilings are exchanged here and in [9] from those in [7] and [8].

(4) Here and in [9], the set of colors is denoted by $\mathbb{A}$ instead of $\Sigma$. Colors are denoted by $a,\ a',\ a_i$ (etc.) instead of $\sigma,\ \sigma',\ \sigma_i$ (etc.).

(5) Here and in [9], the weighted substitution is denoted by $(\sigma, \tau)$ instead of $(\varphi, \eta)$.

(6) Here and in [9], admissible tiles are denoted by $R,\ R',\ R_i,\ R^i$ (etc.) instead of $S,\ S',\ S_i,\ S^i$ (etc.).

(7) Here and in [9], the terminology "primitive" for substitutions is used instead of "mixing" in [7] and [8].

## 3. Weighted substitutions

A *substitution* $\sigma$ on a set $\mathbb{A}$ is a mapping $\mathbb{A} \to \mathbb{A}^+$, where $\mathbb{A}^+ = \bigcup_{\ell=1}^{\infty} \mathbb{A}^\ell$. For $\xi \in \mathbb{A}^+$, we denote $|\xi| := \ell$ if $\xi \in \mathbb{A}^\ell$, and $\xi$ with $|\xi| = \ell$ is usually denoted by $\xi_0\xi_1 \cdots \xi_{\ell-1}$ with $\xi_i \in \mathbb{A}$. We can extend $\sigma$ to be a homomorphism $\mathbb{A}^+ \to \mathbb{A}^+$ as follows:

$$\sigma(\xi) := \sigma(\xi_0)\sigma(\xi_1) \cdots \sigma(\xi_{\ell-1}),$$

where $\xi \in \mathbb{A}^\ell$ and the right-hand side is the concatenations of $\sigma(\xi_i)$'s. We can define $\sigma^2, \sigma^3, \ldots$ as the compositions of $\sigma : \mathbb{A}^+ \to \mathbb{A}^+$.

A *weighted substitution* $(\sigma, \tau)$ on $\mathbb{A}$ is a mapping $\mathbb{A} \to \mathbb{A}^+ \times (0,1)^+$ such that $|\sigma(a)| = |\tau(a)|$ and $\sum_{i<|\tau(a)|} \tau(a)_i = 1$ for any $a \in \mathbb{A}$. Note that $\sigma$ is a substitution on $\mathbb{A}$. We define $\tau^n : \mathbb{A} \to (0,1)^+$ $(n = 2, 3, \ldots)$ (depending on $\sigma$) inductively by

$$\tau^n(a)_k = \tau(a)_i \tau^{n-1}(\sigma(a)_i)_j$$

for any $a \in \mathbb{A}$ and $i, j, k$ with

$$0 \leq i < |\sigma(a)|,\ 0 \leq j < |\sigma^{n-1}(\sigma(a)_i)|,\ k = \sum_{h<i} |\sigma^{n-1}(\sigma(a)_h)| + j.$$

Then, $(\sigma^n, \tau^n)$ is also a weighted substitution for $n = 2, 3, \ldots$.

A substitution $\sigma$ on $\mathbb{A}$ is called *primitive* if there exists a positive integer $n$ such that for any $a, a' \in \mathbb{A}$, $\sigma^n(a)_i = a'$ holds for some $i$ with $0 \le i < |\sigma^n(a)|$.

For a weighted substitution $(\sigma, \tau)$ on $\mathbb{A}$, we always assume that

$$\textit{the substitution } \sigma \textit{ is primitive.} \tag{4}$$

We define the *base set* $B(\sigma, \tau)$ as the closed, multiplicative subgroup of $\mathbb{R}_+$ generated by the set

$$\left\{ \begin{array}{ll} \tau^n(a)_i \ ; & a \in \mathbb{A},\ \ n = 0, 1, \ldots \text{ and } 0 \le i < |\sigma^n(a)| \\ & \qquad\qquad \text{such that } \sigma^n(a)_i = a \end{array} \right\}.$$

**Example 3.1.** Let $\mathbb{A} = \{+, -\}$ and $(\sigma, \tau)$ be a weighted substitution such that

$$\begin{array}{rcl} + & \to & (+, 4/9)(-, 1/9)(+, 4/9) \\ - & \to & (-, 4/9)(+, 1/9)(-, 4/9), \end{array}$$

where we express a weighted substitution $(\sigma, \tau)$ by

$$a \to (\sigma(a)_0, \tau(a)_0)(\sigma(a)_1, \tau(a)_1) \cdots \ (a \in \mathbb{A}).$$

Then, $4/9 \in B(\sigma, \tau)$ since $\sigma(+)_0 = +$ and $\tau(+)_0 = 4/9$. Moreover, $1/81 \in B(\sigma, \tau)$ since $\sigma^2(+)_4 = +$ and $\tau^2(+)_4 = 1/81$. Since $4/9$ and $1/81$ do not have a common multiplicative base, we have $B(\sigma, \tau) = \mathbb{R}_+$. This weighted substitution is discussed in the following sections. The repetition of this weighted substitution starting at $+$ is shown in Figure 3 by colored tiles. The substituted word of a color is represented as the sequence of colors of the connected tiles in below in order from left. The horizontal (additive) sizes of tiles are proportional to the weights and the vertical (multiplicative) sizes are the inverse of the weights.

Let $G := B(\sigma, \tau)$. Then, there exists a function $g : \mathbb{A} \to \mathbb{R}_+$ such that

$$g(\sigma(a)_i)G = g(a)\tau(a)_i G \tag{5}$$

for any $a \in \mathbb{A}$ and $0 \le i < |\sigma(a)|$. Note that if $G = \mathbb{R}_+$, then we can take $g \equiv 1$. In the other case, we can define $g$ by $g(a_0) = 1$ and $g(a) := \tau^n(a_0)_i$ for some $n$ and $i$ such that $\sigma^n(a_0)_i = a$, where $a_0$ is any fixed element in $\mathbb{A}$.

Let $(\sigma, \tau)$ be a weighted substitution satisfying (4). Let $G = B(\sigma, \tau)$. Let $g$ satisfy (5). Let $\Omega(\sigma, \tau, g)'$ be the set of all elements $\omega$ in $\Omega(\mathbb{A})$ such that for any $((x_1, x_2) \times (y_1, y_2), a) \in \omega$, we have

(I) $y_1 \in g(a)G$, and

(II) $(R^i, \sigma(a)_i) \in \omega$ holds for $i = 0, 1, \ldots, |\sigma(a)| - 1$, where

$$\begin{aligned} R^i \ := \ & (x_1 + (x_2 - x_1)\sum_{j=0}^{i-1} \tau(a)_j,\ x_1 + (x_2 - x_1)\sum_{j=0}^{i} \tau(a)_j) \\ & \times (\tau(a)_i y_1,\ y_1). \end{aligned}$$

A vertical line $\gamma := \{x\} \times (-\infty, \infty)$ is called a *separating line* of $\omega \in \Omega(\sigma, \tau, g)'$ if for any $(R, a) \in \omega$, $R \cap \gamma = \emptyset$. Let $\Omega(\sigma, \tau, g)''$ be the set of all $\omega \in \Omega(\sigma, \tau, g)'$ which do not have a separating line and $\Omega(\sigma, \tau, g)$ be the closure of $\Omega(\sigma, \tau, g)''$. Then, the action of $G \times \mathbb{R}$ on $\Omega(\sigma, \tau, g)$ satisfies ($\sharp$1). We usually denote $\Omega(\sigma, \tau, 1)$ simply by $\Omega(\sigma, \tau)$.

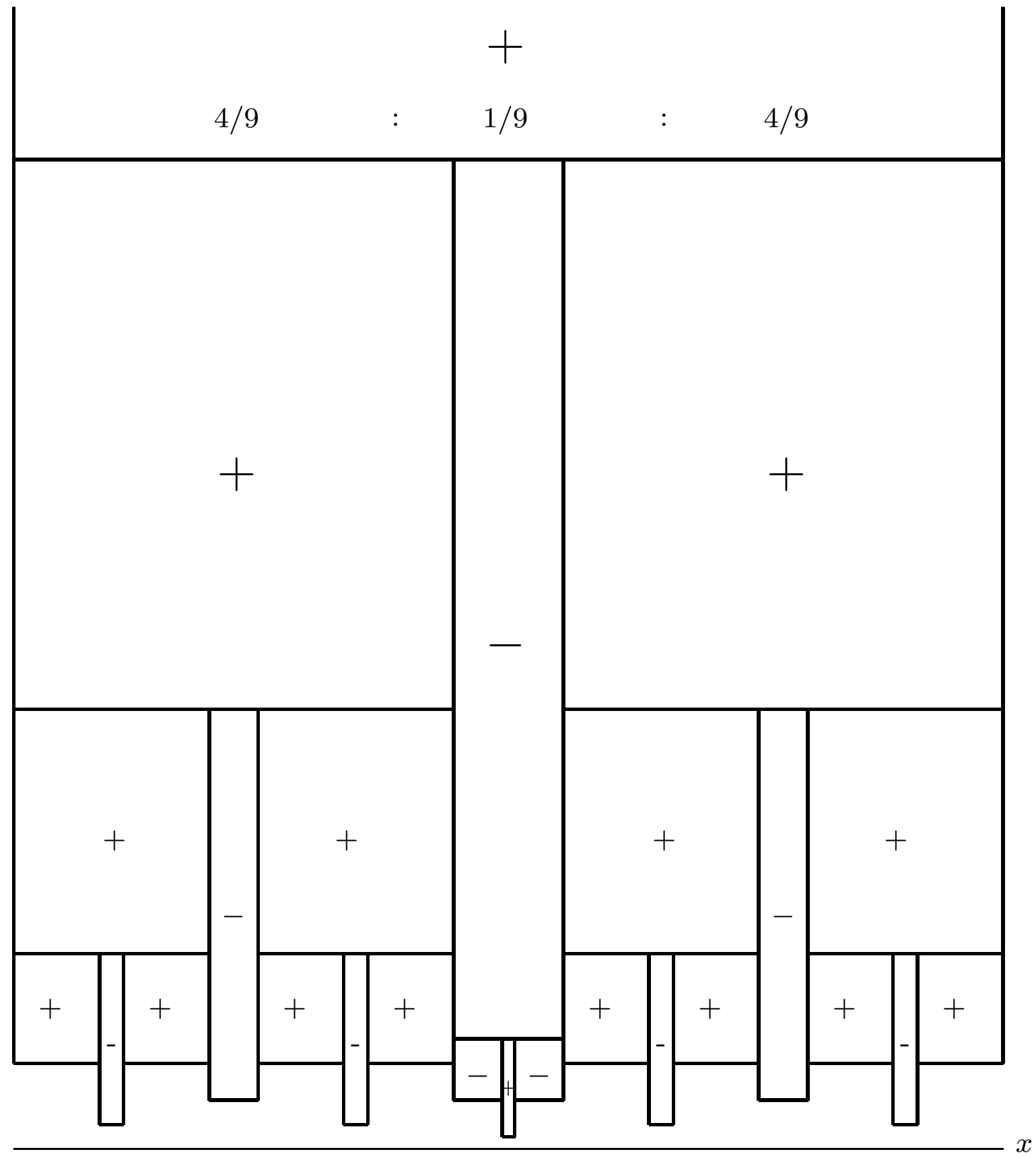

Fig 3. *The weighted substitution in Example 1.*

**Remark 3.2 ([7]).** A nontrivial primitive substitution $\sigma : \mathbb{A} \to \mathbb{A}^+$, where "nontrivial" means $\sum_{a\in\mathbb{A}} |\sigma(a)| \geq 2$, is considered as a weighted substitution in a canonical way. Let

$$M := (\sharp\{0 \leq i < |\sigma(a)|;\ \sigma(a)_i = a'\})_{a,a'\in\mathbb{A}}$$

be the associate matrix. Let $\lambda$ be the maximum eigen-value of $M$ and $\xi := (\xi_a)_{a\in\mathbb{A}}$ be a positive column vector such that $M\xi = \lambda\xi$. Define weight $\tau$ by

$$\tau(a)_i = \frac{\xi_{\sigma(a)_i}}{\lambda\xi_a},$$

which is called the *natural weight* of $\sigma$. Thus, we get a weighted substitution $(\sigma, \tau)$ which admits weight 1. We modify $(\sigma, \tau)$ if necessary in the following way. If there exists $a \in \mathbb{A}$ with $|\sigma(a)| = 1$, so that $a \to (a', 1)$ is a part of $(\sigma, \tau)$, then we replace all the occurrences of $a$ in the right hand side of "$\to$" by $a'$ and remove $a$ from $\mathbb{A}$ together with the rule $a \to (a', 1)$ from $(\sigma, \tau)$. We continue this process until no $a \in \mathbb{A}$ satisfies $|\sigma(a)| = 1$. After that if there exist $a,\ a' \in \mathbb{A}$ such that $(\sigma(a), \tau(a)) = (\sigma(a'), \tau(a'))$, then we identify them.

For example, the 2-*adic expansion substitution* $1 \to 12,\ 2 \to 12$ corresponds to the weighted substitution $1 \to (1, 1/2)(1, 1/2)$. The *Thue-Morse substitution* $1 \to 12,\ 2 \to 21$ corresponds to the weighted substitution $1 \to (1, 1/2)(2, 1/2),\ 2 \to (2, 1/2)(1, 1/2)$. The *Fibonacci substitution* $1 \to 12,\ 2 \to 1$ corresponds to the weighted substitution $1 \to (1, \lambda^{-1})(1, \lambda^{-2})$, where $\lambda = (1 + \sqrt{5})/2$.

The weighted substitution $(\sigma, \tau)$ obtained in this way satisfies that $B(\sigma, \tau) = \{\lambda^n;\ n \in \mathbb{Z}\}$ and that $g$ in (5) can be defined by $g(a) = \xi_a$ $(a \in \mathbb{A})$. Dynamical systems coming from substitutions are discussed by many authors (see [2], for example). Our weighted substitutions are a generalization of them.

Let $(\sigma, \tau)$ be a weighted substitution on $\mathbb{A}$ satisfying (4). Let $g$ satisfy (5). Consider $\Omega(\sigma, \tau, g)$. We call the tile $R^i$ in (II) the $i$-th *child* of the tile $(x_1, x_2) \times (y_1, y_2)$, and $(x_1, x_2) \times (y_1, y_2)$ the *mother* of $R^i$. Note that the vertical size $S(R^i)$ of $R^i$ coincides with the inverse of the weight $\tau(a)_i$. If $R_j$ is a child of $R_{j+1}$ for $j = 0, 1, \ldots, k-1$. Then, the tile $R_0$ is called a $k$-th *descendant* of the tile $R_k$. If $R_0$ is the $i$-th tile among the set of the $k$-th descendants of $R_k$ counting as $0, 1, 2, \cdots$ from left, we call $R_0$ the $(k, i)$-*descendant* of the tile $R_k$. In this case, we also say that $R_k$ is the $k$-th *ancestor* of $R_0$.

**Theorem 3.3.** *The space $\Omega(\sigma, \tau, g)$ is a numeration system with $G = B(\sigma, \tau)$.*

*Proof.* We have already proved ($\sharp$1) and ($\sharp$2) in Theorem 3 of [7]. Here we prove ($\sharp$3). Let $\Omega := \Omega(\sigma, \tau, g)$ and $\mu_\Omega$ be the equilibrium measure. Since $\mu_\Omega$ is the unique invariant probability measure under the additive action, it is also invariant under the multiplicative action.

By Goodman [4], it is sufficient to prove that for any $\lambda \in G$ with $\lambda \neq 1$, the transformation $\omega \mapsto \lambda\omega$ on $\Omega$ has the metrical entropy $|\log \lambda|$ under $\mu_\Omega$, while it has the metrical entropy less than $|\log \lambda|$ under any other $G$-invariant probability measure. □

**Lemma 3.4.** *Let*

$$\begin{aligned} \Sigma &:= \{\omega \in \Omega;\ \omega \textit{ has a separating line}\} \\ \Sigma_0 &:= \{\omega \in \Omega;\ y\textit{-axis is the separating line of } \omega\}. \end{aligned}$$

*Then, we have*

(i) $\Sigma \setminus \Sigma_0$ *is dissipative with respect to the $G$-action, so that $\nu(\Sigma \setminus \Sigma_0) = 0$ for any $G$-invariant probability measure $\nu$ on $\Omega$.*

(ii) *For any $\omega \in \Sigma_0$, $\omega$ restricted to the right quarter plane $(0, \infty) \times (0, \infty)$ and to the left quarter plane $(-\infty, 0) \times (0, \infty)$ are cyclic individually with respect to the $G$-action. Hence, $\overline{G\omega}$ with respect to the $G$-action is either cyclic or conjugate to a 2-dimensional irrational rotation with a multiplicative time parameter.*

(iii) $\Sigma_0$ *is a finite union of minimal and equicontinuous sets with respect to the $G$-action. In fact, there is a mapping from the set of pairs $a \in \mathbb{A}$ and $i$ with $0 \le i < i + 1 < |\sigma(a)|$ onto the set of minimal sets in $\Sigma_0$.*

*Proof.* (i) If the line $x = u$ is the separating line of $\omega \in \Omega$, then $x = \lambda u$ is the separating line of $\lambda\omega$. Hence, $\Sigma \setminus \Sigma_0$ is dissipative.

(ii) Let $\omega \in \Sigma_0$. Denote by $\omega^+$ the restriction of $\omega$ to the right quarter plane $(0, \infty) \times (0, \infty)$, while by $\omega^-$ the restriction of $\omega$ to the left quarter plane $(-\infty, 0) \times (0, \infty)$. Let $(R_i^\pm)_{i \in \mathbb{Z}}$ be the sequence of tiles in $dom(\omega)$ such that $R_i^\pm$ intersects with the upper half lines of $x = \pm 0$, and $R_i^\pm$ is a child of $R_{i+1}^\pm$ for any $i \in \mathbb{Z}$ ($\pm$ respectively). Let $a_i^\pm := \omega(R_i^\pm)$ be the colors of $R_i^\pm$ ($\pm$ respectively). Define mappings $\sigma_\pm$ from $\mathbb{A}$ to $\mathbb{A}$ by $\sigma_+(a) = \sigma(a)_0$ and $\sigma_-(a) = \sigma(a)_{|\sigma(a)|-1}$. Since $\sigma_\pm(a_i^\pm) = a_{i-1}^\pm$ $(i \in \mathbb{Z})$ ($\pm$ respectively), the sequence $(a_i^\pm)_{i \in \mathbb{Z}}$ is periodic, which also implies that the vertical sizes $S(R_i^{\pm})$ of $R_i^{\pm}$, which coincide with the inverses of the weights $\tau(a_{i+1})_{\pm}$, are also periodic in $i \in \mathbb{Z}$ with the period, say $r^{\pm}$ which is the minimum period of $(a_i^{\pm})_{i\in\mathbb{Z}}$ ($\pm$ respectively). Then, $\lambda^+ := \tau^{r^+}(a_0^+)_0^{-1}$ is the minimum (multiplicative) cycle of $\omega^+$, while $\lambda^- := \tau^{r^-}(a_0^-)^{-1}_{|\sigma^{r^-}(a_0^-)|-1}$ is the minimum (multiplicative) cycle of $\omega^-$, that is, $\lambda\omega^{\pm} = \omega^{\pm}$ holds for $\lambda = \lambda^{\pm}$ and $\lambda^{\pm}$ is the minimum among $\lambda > 1$ with this property ($\pm$ respectively).

Therefore, $\omega$ is cyclic with respect to the $G$-action if $\lambda^+$ and $\lambda^-$ have a common multiplicative base. In this case, the minimum cycle of $\omega$ is the minimum positive number $x$ such that $x = (\lambda^+)^n = (\lambda^-)^m$ holds for some positive integers $n, m$. Otherwise, the $G$-action to $\overline{G\omega}$ is conjugate to an 2-dimensional irrational rotation with a multiplicative time parameter.

(iii) We use the notations in the proof of $(ii)$. Take any pair $(a, i)$ with $a \in \mathbb{A}$ and $0 \le i < i + 1 < |\sigma(a)|$. Take any $\omega' \in \Omega$ having a tile $R \in dom(\omega')$ with $\omega'(R) = a$ such that the $y$-axis passes in between the $i$-th child of $R$ and the $i + 1$-th child of $R$. Let $\psi(a, i)$ be the set of limit points of $\lambda\omega'$ as $\lambda \in G$ tends to $\infty$. Note that this does not depend on the choice of $\omega'$. Then, $\psi(a, i)$ is a closed $G$-invariant subset of $\Sigma_0$. Moreover, since the sequence $(\sigma_-^n(\sigma(a)_i), \sigma_+^n(\sigma(a)_{i+1}))_{n=0,1,2,\ldots}$ enter into a cycle after some time, $\psi(a, i)$ is minimal and equicontinuous with respect to the $G$-action.

To prove that the mapping $\psi$ is onto, take any $\omega \in \Sigma_0$. There exists $\omega_n \in \Omega(\sigma, \tau, g)''$ which converges to $\omega$ as $n \to \infty$. We may assume that there exists a pair $(a, i)$ such that for any $n = 1, 2, \cdots$, there exists $R \in dom(\omega_n)$ with $a = \omega_n(R)$ such that the $y$-axis separetes the $i$-th child of $R$ and the $i+1$-th child of $R$. Then, $\omega \in \psi(a, i)$, which proves that $\psi$ is a mapping from the set of pairs $(a, i)$ with $a \in \mathbb{A}$ and $0 \le i < i + 1 < |\sigma(a)|$ onto the set of minimal sets in $\Sigma_0$ with respect to the $G$-action. □

**Example 3.5.** Let $p$ with $0 < p < 1$ satisfy that $\log p/\log(1 - p)$ is irrational. Let $(\sigma, \tau)$ be a weighted substitution on $\mathbb{A} = \{1\}$ such that $1 \to (1, p)(1, 1 - p)$. Then, $B(\sigma, \tau) = \mathbb{R}_+$ holds. Let $\Omega = \Omega(\sigma, \tau)$. In this case, elements in $\Sigma_0$ are not periodic, but almost periodic as shown in Figure 4. Then, the dynamical system $(\Sigma_0,\ \lambda\ (\lambda \in \mathbb{R}_+))$ is isomorphic to $((\mathbb{R}/\mathbb{Z})^2,\ T_\lambda\ (\lambda \in \mathbb{R}_+))$ with

$$T_\lambda(x, y) = (x + \log\lambda/\log(1/p),\ y + \log\lambda/\log(1/(1 - p))).$$

**Lemma 3.6.** *It holds that* $h_{\mu_\Omega}(\lambda) = |\log\lambda|$ *for any* $\lambda \in G$. *Let* $\lambda \ne 1$ *and* $\nu$ *be any other* $\lambda$*-invariant probability measure on* $\Omega$*, then* $h_\nu(\lambda) < |\log\lambda|$.

*Proof.* To prove the lemma, it is sufficient to prove the statements for $\lambda > 1$. Take any $G$-invariant probability measure $\nu$ on $\Omega$ which attains the topological entropy of the multiplication by $\lambda_1 \in G$ with $\lambda_1 > 1$, that is, $h_\nu(\lambda_1) = \log\lambda_1$. We assume also that the $G$-action to $\Omega$ is ergodic with respect to $\nu$. Then by Lemma 3.4, either $\nu(\Sigma_0) = 1$ or $\nu(\Omega \setminus \Sigma) = 1$. In the former case, $h_\nu(\lambda) = 0$ holds for any $\lambda \in G$ since the $G$-action on $\Sigma_0$ is equicontinuous by Lemma 3.4, which contradicts with the assumption. Thus, we have $\nu(\Omega \setminus \Sigma) = 1$.

For $\omega \in \Omega$, let $R_0(\omega) \in dom(\omega)$ be such that $R_0(\omega) = (x_1, x_2) \times (y_1, y_2)$ with $x_1 \le 0 < x_2$ and $y_1 \le 1 < y_2$. Take $a_0 \in \mathbb{A}$ such that

$$\nu(\{\omega \in \Omega;\ \omega(R_0(\omega)) = a_0\}) > 0.$$

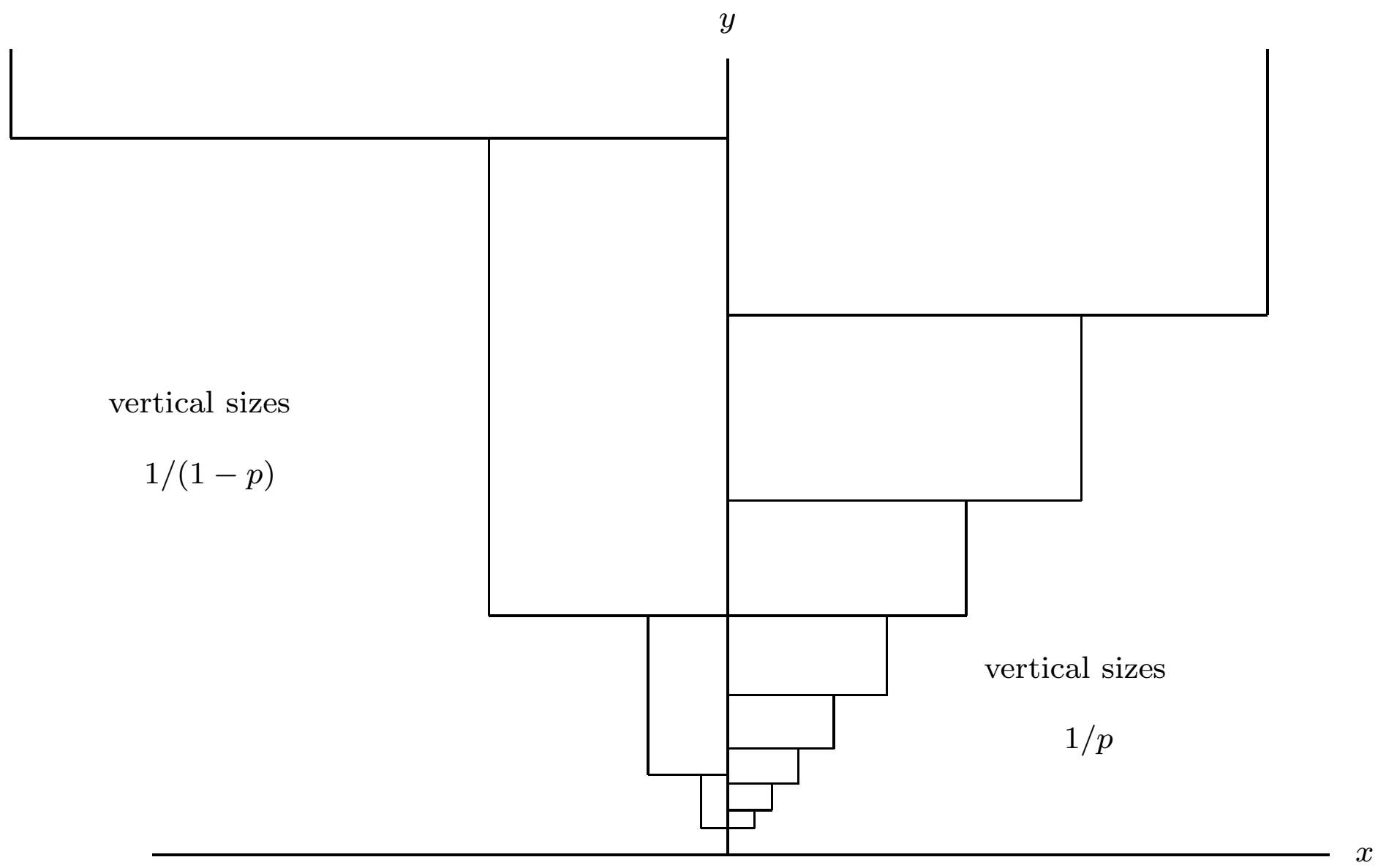

FIG 4. *An element in $\Sigma_0$ in Example 2.*

Take $b_0 := \max\{b \le 1;\ b \in g(a_0)G\}$ (see (5)). Let

$$\begin{aligned}\Omega_1 &:= \{\ \omega \in \Omega;\ \text{the set } \{\lambda \in G;\ \lambda\omega(R_0(\lambda\omega)) = a_0\} \\ &\qquad\qquad \text{is unbounded at } 0 \text{ and } \infty \text{ simultaneously }\} \\ \Omega_0 &:= \{\ \omega \in \Omega_1;\ R_0(\omega) = (x_1, x_2) \times (y_1, y_2) \\ &\qquad\qquad\qquad \text{with } y_1 = b_0 \text{ and } \omega(R_0(\omega)) = a_0\ \}.\end{aligned}$$

For $\omega \in \Omega_0$, let $\lambda_0(\omega)$ be the smallest $\lambda \in G$ with $\lambda > 1$ such that $\lambda\omega \in \Omega_0$. Define a mapping $\Lambda : \Omega_0 \to \Omega_0$ by $\Lambda(\omega) := \lambda_0(\omega)\omega$.

For $k = 0, 1, 2, \cdots$ and $i = 0, 1, \cdots, |\sigma^k(a_0)| - 1$, let

$$\begin{aligned}P(k,i) := \{\ \omega \in \Omega_0;\ &\lambda_0(\omega)^{-1} R_0(\lambda_0(\omega)\omega) \\ &\text{is the } (k,i)\text{-descendant of } R_0(\omega)\}\end{aligned}$$

(see Figure 5) and let

$$\mathcal{P} := \{\ P(k,i);\ k = 1, 2, \cdots,\ 0 \le i < |\sigma^k(a_0)|\ \}$$

be a measurable partition of $\Omega_0$. Note that $\lambda_0(\omega) = \tau^k(a_0)_i^{-1}$ if $\omega \in P(k,i)$.

Since $\nu(\Omega_1) = 1$ by the ergodicity and

$$\Omega_1 = \bigcup_{P(k,i)\in\mathcal{P}} \ \bigcup_{\substack{1\le\lambda<\tau^k(a_0)_i^{-1} \\ \lambda\in G}} \lambda P(k,i),$$

there exists a unique $\Lambda$-invariant probability measure $\nu_0$ on $\Omega_0$ such that for any Borel set $B \subset \Omega$, we have

$$\nu(B) = C(\nu)^{-1} \sum_{P(k,i)\in\mathcal{P}} \int_{b_0}^{b_0\tau^k(a_0)_i^{-1}} \nu_0(\lambda^{-1}B \cap P(k,i))d\lambda/\lambda$$

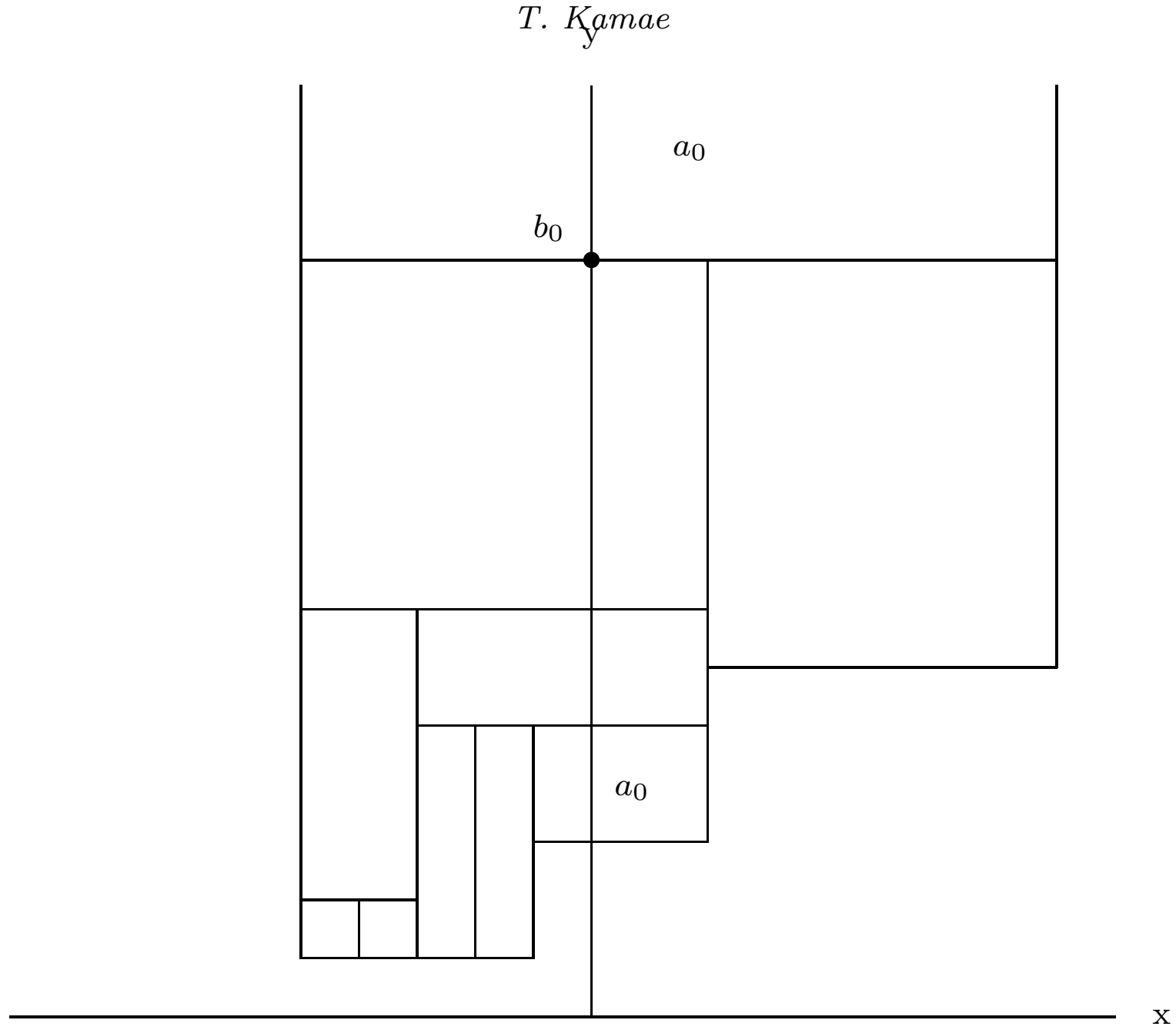

FIG 5. $\omega \in P(3,4)$ *with* $\lambda_0(\omega) = 13/3$

with

$$C(\nu) := \sum_{P(k,i)\in\mathcal{P}} -\log \tau^k(a_0)_i \ \nu_0(P(k,i)) < \infty \tag{6}$$

if $G = \mathbb{R}_+$ and

$$\nu(B) = C(\nu)^{-1} \sum_{P(k,i)\in\mathcal{P}} \sum_{\substack{\lambda\in G \\ b_0 \le \lambda < b_0 \tau^k (a_0)_i^{-1}}} \nu_0(\lambda^{-1}B \cap P(k,i))$$

with

$$C(\nu) := \sum_{P(k,i)\in\mathcal{P}} (-\log \tau^k(a_0)_i / \log \beta) \ \nu_0(P(k,i)) < \infty \tag{7}$$

if $G = \{\beta^n;\ n \in \mathbb{Z}\}$ with $\beta > 1$.

Since

$$\sum_{P(k,i)\in\mathcal{P}} \tau^k(a_0)_i = 1 \text{ and } \sum_{P(k,i)\in\mathcal{P}} \nu_0(P(k,i)) = 1,$$

we have

$$\begin{aligned} H_{\nu_0}(\mathcal{P}) \ &:= \ - \sum_{P(k,i)\in\mathcal{P}} \log \nu_0(P(k,i)) \cdot \nu_0(P(k,i)) \\ &\le \ - \sum_{P(k,i)\in\mathcal{P}} \log \tau^k(a_0)_i \cdot \nu_0(P(k,i)) \end{aligned} \tag{8}$$

by the convexity of $-\log x$. The equality in (8) holds if and only if

$$\nu_0(P(k,i)) = \tau^k(a_0)_i \quad (\forall P(k,i) \in \mathcal{P}). \tag{9}$$

By (6), (7), and (8), we have

$$H_{\nu_0}(\mathcal{P}) = -\sum_{P(k,i)\in\mathcal{P}} \log \nu_0(P(k,i)) \cdot \nu_0(P(k,i)) < \infty.$$

For any $\omega,\ \omega' \in \Omega_0$ such that $\Lambda^k(\omega)$ and $\Lambda^k(\omega')$ belong to the same element in $\mathcal{P}$ for $k = 0, 1, 2, \cdots$, the horizontal position of $R_0(\omega)$, say $(x_1, x_2)$, coincides with that of $R_0(\omega')$. Therefore, $\omega$ and $\omega'$ restricted to $(x_1, x_2) \times (0, b_0)$ coincide. In the same way, if $\Lambda^k(\omega)$ and $\Lambda^k(\omega')$ belong to the same element in $\mathcal{P}$ for any $k \in \mathbb{Z}$, then $R_0 := R_0(\omega) = R_0(\omega')$ holds and all the ancestors of $R_0$ in $\omega$ and $\omega'$ coincide as well as their colors. Therefore, $\omega$ and $\omega'$ restricted to the region covered by the ancestors of $R_0$ coincide. Hence, if $\omega$ or $\omega'$ does not have the separating lines, then $\omega = \omega'$ holds.

Since $\nu(\Sigma) = 0$, we have $\nu_0(\Sigma \cap \Omega_0) = 0$. Hence, the above argument implies that $\mathcal{P}$ is a generator of the system $(\Omega_0, \nu, \Lambda)$. Thus, $h_{\nu_0}(\Lambda) = h_{\nu_0}(\Lambda, \mathcal{P})$. It follows from (8) that

$$\begin{aligned} h_{\nu_0}(\Lambda) &= h_{\nu_0}(\Lambda, \mathcal{P}) \\ &\le H_{\nu_0}(\mathcal{P}) \\ &\le -\sum_{P(k,i)\in\mathcal{P}} \log \tau^k(a_0)_i \cdot \nu_0(P(k,i)). \end{aligned} \tag{10}$$

The equality in the above that

$$h_{\nu_0}(\Lambda) = -\sum_{P(k,i)\in\mathcal{P}} \log \tau^k(a_0)_i \cdot \nu_0(P(k,i))$$

holds if and only if $(\Lambda^n \mathcal{P})_{n\in\mathbb{Z}}$ is an independent sequence with respect to $\nu_0$ satisfying (9).

Since

$$\begin{aligned} h_\nu(\lambda_1)/\log\lambda_1 &= \frac{h_{\nu_0}(\Lambda)}{\int_{\Omega_0} \lambda_0(\omega) d\nu_0(\omega)} \\ &= \frac{h_{\nu_0}(\Lambda)}{-\sum_{P(k,i)\in\mathcal{P}} \log \tau^k(a_0)_i \cdot \nu_0(P(k,i))}\ , \end{aligned}$$

$h_\nu(\lambda_1) \le \log\lambda_1$ follows from (10), while the equality holds if and only if $(\Lambda^n \mathcal{P})_{n\in\mathbb{Z}}$ is an independent sequence with respect to $\nu_0$ satisfying (9). Let this probability measure be $\mu$. Then, it is not difficult to prove that $\mu$ is invariant under the additive action. Hence, the uniqueness of such measure ([7]) proves $\mu = \mu_\Omega$, which completes the proof of Lemma 3.6 and Theorem 3.3. □

The following Theorem 3.7 follows from a known result about the spectrum of unitary operators corresponding to the affine action (Lemma 11.6 of [13], for example).

**Theorem 3.7 ([10]).** *Let $\Omega$ be a numeration system with $G = \mathbb{R}_+$, that is, with the multiplicative $\mathbb{R}_+$-action. Then, the additive action on the probability space $\Omega$ with respect to $\mu_\Omega$ has a pure Lebesgue spectrum.*

## 4. ζ-function

Here, we listed only the results on the $\zeta$-functions. For the proof, refer [10].

Let $\Omega := \Omega(\sigma, \tau, g)$ satisfying (4) and (5). For $\alpha \in \mathbb{C}$, we define the associated matrices on the suffix set $\mathbb{A} \times \mathbb{A}$ as follows:

$$M_\alpha := \left( \sum_{i;\sigma(a)_i = a'} \tau(a)_i^\alpha \right)_{a,a' \in \mathbb{A}} \tag{11}$$

$$M_{\alpha,+} := \left( 1_{\sigma(a)_0 = a'} \ \tau(a)_0^\alpha \right)_{a,a' \in \mathbb{A}}$$

$$M_{\alpha,-} := \left( 1_{\sigma(a)_{|\sigma(a)|-1} = a'} \ \tau(a)_{|\sigma(a)|-1}^\alpha \right)_{a,a' \in \mathbb{A}}$$

Let $\Theta$ be the set of *closed orbits* of $\Omega$ with respect to the action of $G$. That is, $\Theta$ is the family of subsets $\xi$ of $\Omega$ such that $\xi = G\omega$ for some $\omega \in \Omega$ with $\lambda\omega = \omega$ for some $\lambda \in G$ with $\lambda > 1$. We call $\lambda$ as above a *multiplicative cycle* of $\xi$. The minimum multiplicative cycle of $\xi$ is denoted by $c(\xi)$. Note that $c(\xi)$ exists since $\lambda\omega \neq \omega$ for any $\omega \in \Omega$ and $\lambda \in G$ with $1 < \lambda < \min\{\tau(a)_i^{-1};\ a \in \mathbb{A},\ 0 \le i < |\tau(a)|\}$.

We say that $\xi \in \Theta$ has a *separating line* if $\omega \in \xi$ has a separating line. Note that in this case, the separating line is necessarily the $y$-axis and is in common among $\omega \in \xi$. Denote by $\Theta_0$ the set of $\xi \in \Theta$ with the separating line.

Define the $\zeta$-function of $G$-action to $\Omega$ by

$$\zeta_\Omega(\alpha) := \prod_{\xi \in \Theta} (1 - c(\xi)^{-\alpha})^{-1}, \tag{12}$$

where the infinite product converges for any $\alpha \in \mathbb{C}$ with $\mathcal{R}(\alpha) > 1$. It is extended to the whole complex plane by the analytic extension.

**Theorem 4.1.** *We have*

$$\zeta_\Omega(\alpha) = \frac{\det(I - M_{\alpha,+}) \det(I - M_{\alpha,-})}{\det(I - M_\alpha)} \ \zeta_{\Sigma_0}(\alpha),$$

*where*

$$\zeta_{\Sigma_0}(\alpha) := \prod_{\xi \in \Theta_0} (1 - c(\xi)^{-\alpha})^{-1}$$

*is a finite product with respect to* $\xi \in \Theta_0$.

**Theorem 4.2.** (i) $\zeta_\Omega(\alpha) \neq 0$ *if* $\mathcal{R}(\alpha) \neq 0$.

(ii) *In the region* $\mathcal{R}(\alpha) \neq 0$, $\alpha$ *is a pole of* $\zeta_\Omega(\alpha)$ *with multiplicity* $k$ *if and only if it is a zero of* $\det(I - M_\alpha)$ *with multiplicity* $k$ *for any* $k = 1, 2, \ldots$.

(iii) 1 *is a simple pole of* $\zeta_\Omega(\alpha)$.

**Theorem 4.3.** *For* $\Omega = \Omega(\sigma, \eta, g)$, *if* $B(\sigma, \tau) = \{\lambda^n;\ n \in \mathbb{Z}\}$ *with* $\lambda > 1$, *then there exist polynomials* $p, q \in \mathbb{Z}[z]$ *such that* $\zeta_\Omega(\alpha) = p(\lambda^\alpha)/q(\lambda^\alpha)$. *Conversely, if* $\zeta_\Omega(\alpha) = p(\lambda^\alpha)/q(\lambda^\alpha)$ *holds for some polynomials* $p, q \in \mathbb{Z}[z]$ *and* $\lambda > 1$, *then* $B(\sigma, \tau) = \{\lambda^{kn};\ n \in \mathbb{Z}\}$ *for some positive integer* $k$.

**Theorem 4.4.** *If* $B(\sigma, \tau) = \{\lambda^n;\ n \in \mathbb{Z}\}$, *then* $\lambda$ *is an algebraic number.*

## Acknowledgment

The author thanks his old friend, Prof. Mike Keane for his useful discussions and encouragements to develop this research for more than 10 years.